\documentclass[12pt]{article}

\usepackage{amsmath,amssymb,amsthm,amscd}
\usepackage[colorlinks, linkcolor=blue, anchorcolor=gree, citecolor=red]{hyperref}
\usepackage[numbers,sort&compress]{natbib}

\begin{document}

\newcommand{\ts}{\,}
\newcommand{\tss}{\hspace{1pt}}
\newcommand{\Mat}{{\rm{Mat}}}
\newcommand{\CC}{\mathbb{C}}
\newcommand{\Sym}{\mathfrak S}

\newtheorem{thm}{Theorem}[section]
\newtheorem{pro}[thm]{Proposition}
\newtheorem{lem}[thm]{Lemma}
\newtheorem{cor}[thm]{Corollary}
\theoremstyle{definition}
\newtheorem{ex}[thm]{Example}
\newtheorem{remark}[thm]{Remark}
\newcommand{\bth}{\begin{thm}}
\renewcommand{\eth}{\end{thm}}
\newcommand{\bex}{\begin{examp}}
\newcommand{\eex}{\end{examp}}
\newcommand{\bre}{\begin{remark}}
\newcommand{\ere}{\end{remark}}

\newcommand{\bal}{\begin{aligned}}
\newcommand{\eal}{\end{aligned}}
\newcommand{\beq}{\begin{equation}}
\newcommand{\eeq}{\end{equation}}
\newcommand{\ben}{\begin{equation*}}
\newcommand{\een}{\end{equation*}}
\renewcommand{\thefootnote}{}

\newcommand{\bpf}{\begin{proof}}
\newcommand{\epf}{\end{proof}}

\def\beql#1{\begin{equation}\label{#1}}
\title{\Large\bf  On finite groups all of whose cubic Cayley graphs are integral}

\author{{\sc Xuanlong Ma and Kaishun Wang}\\[15pt]
{\small Sch. Math. Sci. {\rm \&} Lab. Math. Com. Sys.,}\\
{\small Beijing Normal University, 100875, Beijing, China.}\\
}
\date{} 

\maketitle

\begin{abstract}
For any positive integer $k$, let  $\mathcal{G}_k$ denote  the set of finite groups $G$  such that all Cayley graphs ${\rm Cay}(G,S)$ are integral whenever $|S|\le k$. Est\'elyi and Kov\'acs \cite{EK14} classified  $\mathcal{G}_k$ for each $k\ge 4$. In this paper, we characterize the finite groups  each of whose cubic Cayley graphs is integral. Moreover, the class $\mathcal{G}_3$ is characterized. As an application, the classification of $\mathcal{G}_k$ is obtained again, where $k\ge 4$.
\end{abstract}



\footnote{E-mail addresses: xuanlma@mail.bnu.edu.cn (X. Ma), wangks@bnu.edu.cn (K. Wang).}

\section{Introduction}

 A graph  is  {\em integral} if  all its eigenvalues are integers.  Harary and Schwenk   \cite{HS74} introduced  integral graphs, and proposed the problem of classifying integral graphs.
Since then   classifications of some special integral graphs have received considerable attention, see \cite{AJ13,BC76,BCL99,BSZ04,CKS11,O09}. For more information, see the two surveys \cite{AABS,BCRSS}.

Let $G$ be a finite group. A subset $S$ of $G$ is called {\em symmetric} if $S^{-1}=S$.
If $S$ is a symmetric subset of $G$ and does not contain the identity, then the {\em Cayley graph}
${\rm Cay}(G,S)$ is the graph with vertex set $G$ and edge set $\{\{g,sg\}: g\in G, s\in S\}$.
Abdollahi and Vatandoost \cite{AV09} listed some infinite families of   integral Cayley graphs, and classified  connected cubic integral Cayley graphs.

A finite group is a {\em Cayley integral group} if each of its
  Cayley graphs  is integral. Klotz and Sander \cite{KS10} introduced this concept and determined all abelian Cayley integral groups. The nonabelian case was handled  by Abdollahi and Jazaeri \cite{AJ14}, and independently by Ahmady et al. \cite{ABM14}.

\begin{thm}\label{allcayig}{\rm (\cite[Theorem 13]{KS10}, \cite[Theorem 1.1]{AJ14}, \cite[Theorem 4.2.]{ABM14})} All finite Cayley integral groups are
$$
\mathbb{Z}_2^m\times \mathbb{Z}_3^n, ~\mathbb{Z}_2^m\times \mathbb{Z}_4^n, ~S_3, ~Q_8\times \mathbb{Z}_2^n, ~{\rm Dic}(\mathbb{Z}_6),
$$
where $m,n \ge 0$, $Q_8$ is the quaternion group of order $8$ and ${\rm Dic}(\mathbb{Z}_6)$ is the generalized dicyclic group of order $12$.
\end{thm}

Very recently, Est${\rm \acute{e}}$lyi and Kov${\rm \acute{a}}$cs \cite{EK14}
generalized this class of groups by introducing the
class $\mathcal{G}_k$ of finite groups $G$  such that all Cayley graphs ${\rm Cay}(G,S)$ are integral whenever $|S|\le k$, and they classified $\mathcal{G}_k$ for each $k\ge 4$.

\begin{thm}\label{mainth}{\rm (\cite[Theorem 1.3]{EK14})}
Every class $\mathcal{G}_k$ consists of the Cayley integral groups if $k\ge 6$. Moreover, $\mathcal{G}_4$ and $\mathcal{G}_5$ are equal, and consist of the Cayley integral groups and
${\rm Dic}(\mathbb{Z}_3^n\times \mathbb{Z}_6)$, where $n\ge 1$.
\end{thm}

For any positive integer $k\ge 2$, let $\mathcal{A}_k$ denote the set of finite groups any of whose Cayley graphs with valency $k$ is integral. We note that If $G$ is group of odd order, then $G\notin \mathcal{A}_m$ for each odd integer $m$.

In this paper we focus on the study of $\mathcal{A}_3$.
In Section \ref{sec:a3}, we characterize $\mathcal{A}_3$, and show that $\mathcal{G}_3$ consists of $\mathcal{A}_3$ and all finite $3$-groups of exponent $3$. As an application, we give an alternating proof Theorem \ref{mainth} in Section \ref{sec:a4}.

\section{Classes $\mathcal{A}_3$ and $\mathcal{G}_3$}\label{sec:a3}

Denote by $D_n$ the dihedral group of order $n$.
Let $V$ be a vector space over complex field $\mathbb{C}$. A {\em representation}
of group $G$ on $V$ is a group homomorphism $\rho$ from $G$ to $GL(V)$, the group of invertible
linear maps from $V$ to itself.
A subspace $W$ of $V$ is said to be {\em invariant} under $\rho$ provided that
$w^{g^{\rho}}\in W$ for any $g\in G$ and $w\in W$. If $V$ has no nontrivial $\rho$-invariant subspaces,
then $\rho$ is called an {\em irreducible representation} of $G$.

\begin{pro}{\rm (\cite[Theorem 3]{DS81})}\label{cayin}
Given a {\rm Cayley} graph ${\rm Cay}(G,S)$, let $$\{\rho_1,\rho_2,\ldots,\rho_t\}$$ be
the set of all irreducible representations of $G$.
Then $\bigcup_{i=1}^{t}\Omega_i$ is the set of all eigenvalues of ${\rm Cay}(G,S)$, where
$\Omega_i$ is the set of all eigenvalues of the matrix
$$\rho_i(S)=\sum_{s\in S}\rho_i(s).$$
\end{pro}


\begin{pro}{\rm(\cite[Theorem 1.1]{AV09})}\label{cubiccig}
A cubic connected {\rm Cayley} graph ${\rm Cay}(G,S)$ is integral
if and only if $G$ is isomorphic one the following groups:
$$
\mathbb{Z}_2^2, \mathbb{Z}_4, \mathbb{Z}_6, \mathbb{Z}_2^3, \mathbb{Z}_2\times \mathbb{Z}_4,
\mathbb{Z}_2\times \mathbb{Z}_6, S_3, D_{8}, D_{12}, A_4, S_4,
D_{8}\times \mathbb{Z}_3, D_{6}\times \mathbb{Z}_4, A_4\times \mathbb{Z}_2.
$$
\end{pro}

Let $G \in \mathcal{A}_k$ and $H$ be a subgroup of $G$.  If $H$ has a subset $S$ such that
$|S|=k$ and $S=S^{-1}$,
since ${\rm Cay}(G,S)$ is a disjoint union of some ${\rm Cay}(H,S)$, one has that
${\rm Cay}(H,S)$ is integral. In particular, if $K \in \mathcal{G}_k$, then
every subgroup of $K$ belongs to $\mathcal{G}_k$. Denote by $\mathcal{G}$ the set of all finite groups $G$ with  $\{|g|: g\in G\}\subseteq\{1,2,3,4,6\}$.

\begin{lem}\label{a2}
A finite group $G$ belongs to $\mathcal{A}_2$ if and only if $G\in \mathcal{G}$, and $G$ is $D_8$-free and $D_{12}$-free.
\end{lem}
\bpf
Suppose that $G\in \mathcal{A}_2$.
It is well known that the cycle of length $n$ is integral if and only if $n=3,4$ or $6$ (cf. \citep[p. 9]{BH12}). This implies that $G\in \mathcal{G}$. Furthermore, if $a$ and $b$ are two generators of $D_8$ such that $|a|=|b|=2$, then ${\rm Cay}(D_8,\{a,b\})$ is the cycle of length $8$ and so $D_8\notin \mathcal{A}_2$. Similarly, we have $D_{12}\notin\mathcal{G}_2$. It follows that $G$ is $D_8$-free and $D_{12}$-free.

For the converse,
let ${\rm Cay}(G,S)$ with valency $2$.
Then $S=\{x,y\}$ or $S=\{z,z^{-1}\}$, where $x$ and $y$ are two distinct involutions, and $z$ is of order greater than $2$.
For the former, one has that $\langle x,y\rangle\cong \mathbb{Z}_2^2$ or $D_6$, since $\mathbb{Z}_2^2$ and  $D_6$ all are Cayley integral,
it follows that ${\rm Cay}(G,\{x,y\})$ is integral. For the latter, obvious that ${\rm Cay}(G,\{z,z^{-1}\})$ is integral.
Thus, we have $G\in \mathcal{A}_2$.
\epf

\begin{lem}\label{ses}
The alternating group $A_4$ belongs to $\mathcal{A}_3$.
\end{lem}
\bpf
Suppose that ${\rm Cay}(A_4,S)$ is a Cayley graph with valency $3$.

\medskip {\em Case 1.} $S$ consists of three involutions.

It is clear that $\langle S\rangle \cong \mathbb{Z}_2^2$. Thereby, we have that ${\rm Cay}(\langle S\rangle,S)$ is integral and so is ${\rm Cay}(A_4,S)$.

\medskip {\em Case 2.} $S$ consists of an involution and two elements of order $3$.

For $i=1$ or $2$, suppose that
$T_i=\{x_i,y_i,y_i^{-1}\}$, where $x_i,y_i\in A_4$, $|x_i|=2$, $|y_i|=3$, and.
Then $\langle T_i\rangle=A_4$, and it is easy to check that the mapping
$\sigma:  x_1\longmapsto x_2, y_1\longmapsto y_2$ is an automorphism of $A_4$. This means that
$\sigma$ is an isomorphism from ${\rm Cay}(A_4,T_1)$ to ${\rm Cay}(A_4,T_2)$. Consequently, to see $A_4 \in \mathcal{A}_3$, it is sufficient to prove that ${\rm Cay}(A_4,T_1)$ is integral.
Take $T_1=S=\{(1,3)(2,4),(2,4,3),(2,3,4)\}$ and
let $\omega=e^{\frac{2\pi}{3}i}$.
By GAP \cite{G13} all nontrivial irreducible representations of  $A_4$ are given by
$$
\rho_1: (2,4,3)\mapsto \omega, (1,3)(2,4)\mapsto 1; \rho_2: (2,4,3)\mapsto \omega^2, (1,3)(2,4)\mapsto 1;
$$
$$
\rho_3: (2,4,3)\longmapsto \left(
                         \begin{array}{ccc}
                           0 & 0 & 1 \\
                           1 & 0 & 0 \\
                           0 & 1 & 0 \\
                         \end{array}
                       \right),
(1,3)(2,4)\longmapsto \left(
                         \begin{array}{ccc}
                           -1 & 0 & 0 \\
                           0 & 1 & 0 \\
                           0 & 0 & -1 \\
                         \end{array}
                       \right).
$$
By Proposition \ref{cayin}, it is easy to check that ${\rm Cay}(A_4,S)$ is integral.
\epf

\begin{pro}\label{c3lem2}
Let $G\in\mathcal{A}_3$. Then $G$ has a subgroup isomorphic to $S_3$ if and only if $G\cong S_3$.
\end{pro}
\bpf
Assume that $K$ is a subgroup of $G$ isomorphic to $S_3$. Let $\{x,y,z\}$ be the set of all involutions of $K$.

Suppose that $a$ is an involution in $G\setminus \{x,y\}$. Then $\langle x,y,a\rangle$ is nonabelian and
${\rm Cay}(\langle x,y,a\rangle,\{ x,y,a\})$ is a cubic connected integral graph.
Note that $\langle x,y\rangle\cong S_3$.
It follows from Proposition \ref{cubiccig} that $\langle x,y,a\rangle$ is one the following groups:
\begin{equation}\label{f1}
S_3, ~D_{12}, ~S_4, ~D_{6}\times \mathbb{Z}_4.
\end{equation}

We now claim $D_8\notin \mathcal{A}_3$.
Let $D_8=\langle a,b: a^4=b^2=1, bab=a^3\rangle$. Then by GAP \cite{G13},
$D_8$ has a $2$-dimensional irreducible representation $\rho_1$ given by
$$
a\longmapsto \left(
               \begin{array}{cc}
                 0 & -1 \\
                 1 & 0 \\
               \end{array}
             \right), b\longmapsto \left(
               \begin{array}{cc}
                 0 & 1 \\
                 1 & 0 \\
               \end{array}
             \right).
$$
Set $T_1=\{a^2,a^3b,b\}$. Then we have
$$
\sum_{s\in T_1}\rho_1(s)=\left(
                            \begin{array}{cc}
                              0 & 1 \\
                              1 & -2 \\
                            \end{array}
                          \right),
$$
which implies that $-1\pm \sqrt{2}$ is an eigenvalue of ${\rm Cay}(D_8,T_1)$ by Proposition \ref{cayin}. Thus, $D_8\notin \mathcal{A}_3$.

Now we prove that $D_{12}\notin \mathcal{A}_3$. Let $D_{12}=\langle a,b: a^6=b^2=1, bab=a^3\rangle$. Then $D_{12}$ has a $2$-dimensional irreducible representation $\rho_2$ determined by
$$
a^5b\longmapsto \left(
               \begin{array}{cc}
                 0 & \omega \\
                 \omega^2 & 0 \\
               \end{array}
             \right), b\longmapsto \left(
               \begin{array}{cc}
                 0 & -1 \\
                 -1 & 0 \\
               \end{array}
             \right).
$$
Take $T_2=\{a^3,a^5b,b\}$. It follows that
$$
\sum_{s\in T_2}\rho_2(s)=\left(
                            \begin{array}{cc}
                              -1 & \omega-1 \\
                              \omega^2-1 & -1 \\
                            \end{array}
                          \right),
$$
which  has  eigenvalues $-1\pm \sqrt{3}$ and so one has that $D_{12}\notin \mathcal{A}_3$.

Note that $S_4$ contains a subgroup isomorphic to $D_8$, and $D_6\times \mathbb{Z}_4$ has a subgroup isomorphic to $D_{12}$. Since $\langle x,y,a\rangle$  belongs to $\mathcal{A}_3$,
by (\ref{f1}) one has that  $\langle x,y,a\rangle\cong S_3$.
This means that $a=z$.
Thus, $G$ has precisely three distinct involutions and so $K$
is normal in $G$.
Note that if $G$ has an element $g$ of order greater than $2$, then ${\rm Cay}(\langle x,g,g^{-1}\rangle,\{x,g,g^{-1}\})$ is a cubic connected integral graph. By Proposition \ref{cubiccig}, it is easy to see that $G\in \mathcal{G}$.

Suppose that $b$ is an element of $G$ with $|b|=4$. Let $g$ be an arbitrary involution of $G$.
If $\langle g,b\rangle$ is nonanelian, then $\langle g,b\rangle\cong A_4$ or $S_3$, a contradiction since $A_4$ and $S_3$ have no elements of order $4$. Hence, one gets $[g,b]=1$,
where $[g,b]$ is the commutator of $g$ and $b$, that is, $[g,b]=g^{-1}b^{-1}gb$.
It follows that $[b^2,g]=1$, which is impossible. Thus, $G$ has no elements of order $4$.
By a similar argument, $G$ also has no elements of order $6$.
Thereby, it follows that $\{|g|: g\in G\}\subseteq\{1,2,3\}$.

Let $w$ belong to the centralizer of $K$ in $G$. If $|w|=3$ then $G$ has an element
$wx$ with order $6$, a contradiction. Furthermore,
since all  involutions are pairwise noncommutative,  one has $|w|\ne 2$. This means that $C_G(K)=1$. By the $N/C$ Theorem (cf. \citep[Theorem 1.6.13]{Rob}), we obtain that $G$ is isomorphic to a subgroup of the full automorphism group ${\rm Aut}(K)$ of $K$. Note that ${\rm Aut}(K)\cong S_3$.  Thus, we conclude that $G\cong S_3$.

The converse implication is straightforward. \epf

By the definition of $\mathcal{A}_3$, we see that every group of odd order does not belong to  $\mathcal{A}_3$. Now we give a characterization for class $\mathcal{A}_3$.

\begin{thm}\label{thm2} Let $G$ be a finite group. Then $G\in\mathcal{A}_3$ if and only if $G\cong S_3$, or for any involution $x$ and element $y\in G$, $\langle x,y\rangle$ is isomorphic to one of the following groups:
\begin{equation}\label{equ1}
\mathbb{Z}_2, ~\mathbb{Z}_2^2, ~\mathbb{Z}_4, ~\mathbb{Z}_6, ~\mathbb{Z}_2\times \mathbb{Z}_4,
~\mathbb{Z}_2\times \mathbb{Z}_6, ~A_4.
\end{equation}
\end{thm}
\bpf
$``\Rightarrow"$: {\em Case 1}. $|y|=1$.

It is clear that $\langle x,y\rangle\cong \mathbb{Z}_2$.

{\em Case 2}. $|y|=2$.

If
$y=x$ or $[x,y]=1$, then $\langle x,y\rangle\cong \mathbb{Z}_2$ or $\mathbb{Z}_2^2$, as desired.
Note that two distinct involutions generate a dihedral group.
Thus, we may suppose that $\langle x,y\rangle\cong D_{2n}$, where $n\ge 3$.
This implies that ${\rm Cay}(\langle x,y\rangle,\{x,y,z\})$ is integral, where $z$ is an involution of
$\langle x,y\rangle\setminus \{x,y\}$.
Since $D_8,D_{12}\notin \mathcal{A}_3$,   one has that $\langle x,y\rangle\cong S_3$ by Proposition \ref{cubiccig}.

{\em Case 3}. $|y|\ge 3$.

Then
${\rm Cay}(\langle x,y,y^{-1}\rangle,\{x,y,y^{-1}\})$ is a cubic connected integral graph.
Note that $S_4,D_8\times \mathbb{Z}_3,D_6\times \mathbb{Z}_4,A_4\times \mathbb{Z}_2\notin \mathcal{A}_3.$
By Proposition \ref{cubiccig} again, one has that $\langle x,y\rangle\cong \mathbb{Z}_4, \mathbb{Z}_6, \mathbb{Z}_2\times \mathbb{Z}_4 , \mathbb{Z}_2\times \mathbb{Z}_6, S_3$ or $A_4$.

From the above, now the desired result follows from Proposition \ref{c3lem2}.

$``\Leftarrow"$: It is clear that $S_3\in \mathcal{A}_3$. We may assume that $G\ncong S_3$.
Let $S$ be an arbitrary symmetric subset
of $G$ with $|S|=3$ and $1\notin S$.

{\em Case 1}. $S$ consists of three involutions.

Note that  every two elements of $S$ commute. So $\langle S\rangle\cong \mathbb{Z}_2^2$ or $\mathbb{Z}_2^3$.
This means that ${\rm Cay}(G,S)$ is integral.

{\em Case 2}. $S=\{x,y,y^{-1}\}$ with $|x|=2$ and $|y|\ge3$.

If $[x,y]\neq 1$, then $\langle x,y\rangle\cong A_4$, by Lemma \ref{ses} one concludes that
${\rm Cay}(G,S)$ is integral. If $\langle x,y\rangle$ is abelian, then $\langle x,y\rangle$ is Cayley integral and so is ${\rm Cay}(G,S)$. \epf

Note that all $3$-groups of exponent $3$ are contained in $\mathcal{G}_3$, but they are not contain in $\mathcal{A}_3$. Combining  Lemma \ref{a2} and Theorem \ref{thm2}, we obtain a characterization of
$\mathcal{G}_3$.

\begin{cor}\label{g3}
The class $\mathcal{G}_3$ consists of all finite $3$-groups of exponent $3$ and all finite groups in  $\mathcal{A}_3$.
\end{cor}

\begin{cor}\label{nil-g3}
Let $G$ be a nilpotent group.  Then $G\in\mathcal{G}_3$ if
and only if one of the following holds:

(1) $G$ is a $3$-group of exponent $3$;

(2) $G\cong \mathbb{Z}_2^n$ for $n\ge 1$;

(3) $G$ is a $2$-group of exponent $4$, and every involution of $G$ belongs to $Z(G)$, the center of $G$;

(4) $G\cong \mathbb{Z}_2^n\times B$, where $B$ is an  arbitrary $3$-group with exponent $3$, and $n\ge 1$.
\end{cor}

As pointed out in \cite{K14}, $\mathcal{G}_3$ is much wider than the classes $\mathcal{G}_k$ for $k\ge 4$.
We now present some examples belonging to $\mathcal{G}_3$, however, they all do not belong to $\mathcal{G}_4$. Firstly, by Corollary \ref{g3} we see that $A_4\in \mathcal{G}_3$.

\begin{ex}\label{nil2g}
Let $A$ be an abelian group with exponent $4$. Then $Q_8\times A\in \mathcal{G}_3$. It is because that for any involution $t$ of $Q_8\times A$, we have that $t=(-1,1),(1,x)$ or $(-1,x)$, where $x$ is an involution of $A$. Thereby, $t\in Z(Q_8)\times A=Z(Q_8\times A)$. In view of Corollary \ref{nil-g3}, one has that $Q_8\times A\in \mathcal{G}_3$.
\end{ex}

\begin{ex}\label{A4XZ3}
Let $B$ be a group with exponent $3$ and
let $A_4\setminus \{(1)\}=T\cup H$, where $T$ is the set of all involutions and $H$ is the set of all elements of order $3$.
Note that every involution in $A_4\times B$ has the form $(t,1)$ for some $t\in T$. It is easy to see that if $t_0,t\in T$ and $h\in H$, then
$[(t,1),(t_0,1)]=1$, $[(t,1),(t_0,a)]=1$ and $\langle(t,1),(h,1)\rangle\cong A_4$, where $1\ne a\in B$.
Now take $t\in T$ and $h\in H$, one has that
$$\langle(t,1),(h,a): (t,1)^2=(h,a)^3=1, ((h,a)(t,1))^3=1\rangle \cong A_4.$$
Thus, we have that $A_4\times B\in \mathcal{G}_3$ by Corollary \ref{g3}.
\end{ex}

\begin{ex}\label{sl23}
By Corollary \ref{g3}, if $G\in \mathcal{G}$ and every involution of $G$ is central, then $G\in \mathcal{G}_3$. For example, the special linear group $SL(2,3)$, clearly, $\{|g|: g\in SL(2,3)\}=\{1,2,3,4,6\}$ and $SL(2,3)$
has precisely one involution, so $G\in\mathcal{G}_3$. Particularly, $\mathbb{Z}_2^n\times SL(2,3)\in \mathcal{G}_3$ for each $n\ge 2$.
\end{ex}

\section{Proof of Theorem \ref{mainth}}\label{sec:a4}

In this section, by using Corollary \ref{g3}, we give an alternative proof of Theorem \ref{mainth}.
By \cite[Lemma 2.4]{EK14}, for each $n\ge 1$, the groups {\rm Dic}$(\mathbb{Z}_3^n\times \mathbb{Z}_6)$ belongs to $\mathcal{G}_5\setminus \mathcal{G}_6$. In order to prove Theorem \ref{mainth}, we only need
to show that any group in $\mathcal{G}_4$ is a Cayley integral group or {\rm Dic}$(\mathbb{Z}_3^n\times \mathbb{Z}_6)$, where $n$ is a positive integer.

\begin{lem}\label{a41}
Let $G$ be a nilpotent group. If $G\in\mathcal{G}_4$, then $G$ is {\rm Cayley} integral.
\end{lem}
\bpf
Note that $\mathcal{G}_4\subseteq\mathcal{G}_2$. Then, by Lemma \ref{a2} we have that $G\in \mathcal{G}$. It means that $|G|$ has at most two distinct prime divisor $2$ and $3$.

\medskip {\em Case 1.} $G$ is a $3$-group.

If $G$ is abelian, then $G$ is elementary abelian and so $G$ is Cayley integral, as desired.

Suppose that $G$ is nonabelian.
Since $G$ is of exponent $3$, $G$ has two elements $a,b$ such that $|a|=|b|=3$ and $[a,b]\ne 1$. It follows that $\langle a,b\rangle$ is the nonabelian group of order $27$ and exponent $3$, where $\langle a,b\rangle=\langle a,b: a^3=b^3=(ab)^3=(ab^2)^3=1\rangle$.
Take $S=\{a,a^2,b,b^2\}$.  Note that $\langle a,b\rangle$ has a $3$-dimensional irreducible representation given by
$$
a\longmapsto \left(
    \begin{array}{ccc}
      1 & 0 & 0 \\
      0 & \omega & 0 \\
      0 & 0 & \omega^2 \\
    \end{array}
  \right),
b\longmapsto \left(
    \begin{array}{ccc}
      0 & 0 & \omega^2 \\
      \omega & 0 & 0 \\
      0 & 1 & 0 \\
    \end{array}
  \right).
$$
It is easy to check that ${\rm Cay}(\langle a,b\rangle,S)$ is not integral,
which is a contradiction as $\langle a,b\rangle\in\mathcal{G}_4$.

\medskip {\em Case 2.} $G$ is a $2$-group.

It is clear that $G$ is Cayley integral if $G$ is abelian.
Thus, we may
assume that $G$ is nonabelian. Then $G$ is of exponent $4$ and
by Corollary \ref{nil-g3}, one has that every cyclic subgroup of order $2$ is
central in $G$.

Now we show that every cyclic subgroup of order $4$ is also normal in $G$. To see this,
suppose, to the contrary, that there exist two elements $x$ and $y$ in $G$ such that $|x|=4$
and $x^y\notin \langle x\rangle$. Then $|y|=4$ and $[x,y]\ne 1$.
Considering
$$H=\langle x,y: x^4=y^4=[x^2,y]=[x,y^2]=(xy)^4=1\rangle,$$
with the help of GAP \cite{G13}, one concludes that
$H\cong H_0,H_1,H_2$ or $Q_{8}$,
where $$H_{0}=\langle a,b,c: a^4=b^4=c^2=1,[a,b]=c,[c,a]=[c,b]=1\rangle,$$
$$H_1=\langle a_1,b_1: a_1^4=b_1^4=1,a_1^{b_1}=a_1^{-1}\rangle,$$
$$H_2=\langle a_2,b_2,c_2: a_2^4=b_2^2=c_2^2=1,[a_2,b_2]=c_2,[c_2,a_2]=[c_2,b_2]=1\rangle.$$

Now we prove that $H_0$, $H_{1}$ and $H_2$ do not belong to $\mathcal{G}_4$.
Note that $b_2$ is an involution of $H_2$, and $b_2\notin Z(H_2)$. Consequently, we have that $H_2\notin \mathcal{G}_3$ by Corollary \ref{nil-g3}, and hence $H_2\notin \mathcal{G}_4$.
For $H_1$,
set $T_1=\{a_1^2b_1^{-1},b_1a_1^2,a_1^{-1}b_1^{-1},b_1a_1\}$, since $H_1$ has a $2$-dimensional irreducible representation given by
$$
a_1\longmapsto \left(
                 \begin{array}{cc}
                   0 & -1 \\
                   1 & 0 \\
                 \end{array}
               \right),
               b_1\longmapsto \left(
                 \begin{array}{cc}
                   0 & 1 \\
                   1 & 0 \\
                 \end{array}
               \right),
$$
it is easy to check that ${\rm Cay}(H_1,T_1)$ is not integral. For $H_0$, by verifying one see that
$H_{0}=\langle ba^2,a^3b^2\rangle$, take $T_0=\{ba^2,a^2b^3,a^3b^2,b^2a\}$. Note that  $H_0$ has a $2$-dimensional irreducible representation given by
$$
ba^2\longmapsto \left(
                  \begin{array}{cc}
                    0 & 1 \\
                    1 & 0 \\
                  \end{array}
                \right),
                a^3b^2\longmapsto \left(
                  \begin{array}{cc}
                    1 & 0 \\
                    0 & -1 \\
                  \end{array}
                \right).
$$
This implies that ${\rm Cay}(H_{0},T_0)$ is not integral.

Considering the above we see that
$H\cong Q_8$, which is a contradiction since every subgroup of $Q_8$ is normal. This yields that every cyclic subgroup of order $4$ is normal in $G$.

Now note that every cyclic subgroup of $G$ is normal and so is every subgroup of $G$. Consequently, $G$ is isomorphic to a direct product of $Q_8$,
an elementary abelian $2$-group and an abelian group of odd order (cf. \cite{B33}), one gets that
$G\cong Q_8\times \mathbb{Z}_2^m$ for some nonnegative integer $m$. By Theorem \ref{allcayig} one has that $G$ is Cayley integral, as desired.

\medskip {\em Case 3.} $G=P\times Q$, where $P$ and $Q$ are the Sylow $2$- and $3$-subgroups of $G$, respectively.

By Case 1, one has that $Q$ is elementary abelian. Since $G$ has no elements of order $12$, one gets that $P$ is also elementary abelian. It means that $G$ is abelian, and so $G$ is  Cayley integral. \epf

\begin{lem}\label{a42}
Let $G$ be a nonnilpotent group. If $G\in\mathcal{G}_4$, then $G\cong S_3$ or $G\cong {\rm Dic}(\mathbb{Z}_3^n\times \mathbb{Z}_6)$, where $n$ is a nonnegative integer.
\end{lem}
\bpf
We first claim that $A_4\notin \mathcal{G}_4$. Clearly, $A_4=\langle (2,4,3),(1,3)(2,4),(1,2)(3,4)\rangle$. Take $S=\{(2,3,4),(2,4,3),(1,3)(2,4),(1,2)(3,4)\}$. Note that $A_4$ has  a $3$-dimensional irreducible representation determined by
$$
(2,4,3)\longmapsto \left(
                     \begin{array}{ccc}
                       0 & 0 & 1 \\
                       1 & 0 & 0 \\
                       0 & 1 & 0 \\
                     \end{array}
                   \right),
(1,3)(2,4)\longmapsto \left(
                     \begin{array}{ccc}
                       -1 & 0 & 0 \\
                       0 & 1 & 0 \\
                       0 & 0 & -1 \\
                     \end{array}
                   \right),
$$
$$
(1,2)(3,4)\longmapsto \left(
                     \begin{array}{ccc}
                       1 & 0 & 0 \\
                       0 & -1 & 0 \\
                       0 & 0 & -1 \\
                     \end{array}
                   \right).
$$
By verifying one has that ${\rm Cay}(A_4,S)$ is not integral, so our claim is valid.

Note that $G\in \mathcal{G}_3$.
By Corollary \ref{g3} and Proposition \ref{c3lem2}, we may assume that $G\ncong S_3$.
Then again by Corollary \ref{g3} and $A_4\notin \mathcal{G}_4$, every involution of $G$ belongs to $Z(G)$. Since $G$ is not nilpotent, $G$ has elements of order $4$.
Suppose that $b$ and $a$ are two elements of $G$ such that $|b|=4$ and $|a|=3$.
Note that $G$ has no elements of order $12$. Then $[a,b]\ne 1$.
Since $G\in \mathcal{G}_4$, ${\rm Cay}(\langle a,b\rangle,\{a,a^{-1},b,b^{-1}\})$ is a quartic connected integral graph. Note that all quartic connected Cayley integral graphs were obtained in \cite{MW15}. Observe that $a^b=a^{-1}$.

Now we claim that $G$ has a unique involution.
Suppose, to the contrary, that there exists an
involution $u$ in $G$ such that
$u\ne b^2$. Then $G$ has a subgroup
$$
H=\langle a,b,u: a^3=b^4=u^2=1, [a,u]=[b,u]=1, a^b=a^{-1}\rangle\cong (\mathbb{Z}_3\rtimes \mathbb{Z}_4)\times \mathbb{Z}_2.
$$
Let $S=\{b^{-1}u, ub, ba, a^{-1}b^{-1}\}$. Note that there exists a irreducible representation of $H$
given by
$$
a^{-1}\longmapsto \left(
                  \begin{array}{cc}
                    \omega & 0 \\
                    0 & \omega^2 \\
                  \end{array}
                \right),
b\longmapsto \left(
                  \begin{array}{cc}
                    0 & -1 \\
                    -1 & 0 \\
                  \end{array}
                \right),
u\longmapsto \left(
                  \begin{array}{cc}
                    -1 & 0 \\
                    0 & -1 \\
                  \end{array}
                \right).
$$
It is easy to check that ${\rm Cay}(H,S)$ is not integral, contrary to $H\in \mathcal{A}_4$.
This forces that the claim is valid and thereby, $G$ has a Sylow $2$-subgroup $P$ that is isomorphic to $\mathbb{Z}_4$ or $Q_8$; this is because that $P$ has a unique subgroup of order $2$ (cf. \cite[pp. 252, Theorem 9.7.3]{S64}). Write
$$Q_8=\{\pm1, \pm i, \pm j, \pm k: i^2=j^2=k^2=-1, ij=k=-ji\}.$$
Note that $|a|=3$.
If $P= Q_8$ then $a^{ij}=a^{-1}=(a^{-1})^j$ and so $|a^{-1}j|=12$, a contradiction. It follows that
$P=\langle b\rangle$. Note that $G$ has a subgroup $Q$ isomorphic to $\mathbb{Z}_3^m\times \mathbb{Z}_2$ for $m\ge 1$, and for each $x\in Q$ one has $x^b=x^{-1}$.
Thus, $G$ is isomorphic to Dic$(\mathbb{Z}_3^n\times \mathbb{Z}_6)$ for $n\ge 0$.
\epf

Combining Lemma \ref{a41} and Lemma \ref{a42}, we complete the proof of Theorem \ref{mainth}.

\section{Acknowledgement}
The authors are grateful to the referee for a very careful reading of the paper, and
many useful suggestions and comments.
This research is supported by National Natural Science Foundation of China (11271047,
11371204) and the Fundamental Research Funds for the Central University of China.

\end{document}